\documentclass{elsart}
 \usepackage{amssymb} 
 \usepackage{graphicx} 

\def\proofend{$\Box$\vspace{0.4em}}

\def\Aut{\mbox{Aut}}

\def\proofr{{\indent\it P\,r\,o\,o\,f\,.\,~}}

\theoremstyle{plain} 
\newtheorem{thrm}{Theorem}[section]
\newtheorem{prp}[thrm]{Proposition}
\newtheorem{crll}[thrm]{Corollary}
\newtheorem{lmm}[thrm]{Lemma}
\newtheorem{exs}[thrm]{Exercise}
\theoremstyle{definition} 
\newtheorem{rmrk}[thrm]{Remark}
\newtheorem{xmpl}[thrm]{Example}

\newcommand\btheorem{\begin{thrm}}
\newcommand\etheorem{\end{thrm}}

\newcommand\blemma{\begin{lmm}}
\newcommand\elemma{\end{lmm}}

\newcommand\bpro{\begin{prp}}
\newcommand\epro{\end{prp}}

\newcommand\bcorol{\begin{crll}}
\newcommand\ecorol{\end{crll}}

\newcommand\bexs{\begin{exs}}
\newcommand\eexs{\end{exs}}

\newcommand\bnote{\begin{rmrk}}
\newcommand\enote{\end{rmrk}}

\newcommand\bexam{\begin{xmpl}}
\newcommand\eexam{\end{xmpl}}

\journal{Discrete Mathematics}
\begin{document}
\begin{frontmatter}
\title{On diameter perfect constant-weight ternary codes}
\author{
D. S. Krotov
}
\ead{krotov@math.nsc.ru}
\address{Sobolev Institute of Mathematics, Koptyuga, 4,
Novosibirsk, 630090, Russia}
\begin{abstract}
From cosets of binary Hamming codes we construct
diameter perfect constant-weight ternary codes with weight $n-1$
(where $n$ is the code length) and distances $3$ and $5$.
The class of distance $5$ codes has parameters unknown before.
\end{abstract}
\end{frontmatter}

\section{Introduction}
Let $S$ be a finite metric space.
A subset $C$ of $S$ is called a {\it distance $d$ code} if $|C|\geq 2$
and $d$ is the minimal distance between any two different words in $C$.
A nonempty subset $A$ of $S$ is called a {\it diameter $D$ anticode}
if $D$ is the maximal distance between any two words in $A$.
The concept of diameter perfect code \cite{AhlAydKha} is based on a
corollary of the following well known fact.

\blemma \label{cslb} Assume the group of isometries of the metric space $S$ is transitive.
Let $C$ be a code in $S$ with distances from
$\bar d=\{d_1,...,d_k\}$. Further let $L \subset B$ be a maximal code in
$B \subset S$ with distances from $\bar d$. Then one has
\[
|C|/|S| \leq |L|/|B|.
\]
\elemma
\bcorol  Assume the group of isometries of the metric space $S$ is transitive.
Let $C$ be a distance $d$ code,
  $A$ be a diameter $D$ anticode in $S$ and $D<d$. Then
\begin{equation}\label{CA}
  |C| \cdot |A| \leq |S|.
\end{equation}
\ecorol
\proofr
Take $B=A$ and $\bar d=\{d'\,|\,d'\geq d\}$ in Lemma \ref{cslb}.
Then $|L|=1$ and the statement follows.
\proofend
\bnote As was found by Delsarte {\rm\cite{Delsarte:1973}}, the bound {\rm(\ref{CA})}
also holds for
metric spaces generated by distance regular graphs
{\rm (}even if the group of isometries is not transitive{\rm )}.
\enote

If (\ref{CA}) holds with equality, then the
code $C$ is called a {\it diameter perfect distance $d$ code} or {\it $D$-diameter perfect
code}.
It is obvious that in this case $A$ is a maximal (by cardinality)
diameter $D$ anticode.

In this paper we consider the
space of ternary $n$-words of weight $n-1$ with Hamming metric.
Recall that the Hamming distance $d_H$ between two $n$-words is the number of
positions in which they differ.
For convenience,  we replace the symbols of
the ordinary ternary alphabet $\{0,1,2\}$ using the substitutions $0\to *$, $1\to
0$, $2\to 1$.
So, we get the space $X^n$ defined as the set of $n$-words over the alphabet $\{0,1,*\}$
with exactly one symbol $*$. Note that $|X^n|=n2^{n-1}$.

In Section \ref{3} we construct diameter perfect ternary constant-weight codes
with distance $3$. Codes with such parameters was known before, see
\cite{Svanstr:PhD,Svanstr:1999,vLinTol:1999,
Kro:2001:Ternary}.

Section \ref{4} contains some notes on diameter perfect ternary constant-weight codes
with distance $4$.

In Section \ref{5} we construct a class of
diameter perfect distance $5$ code in $X^n$ for $n=2^m$ where $m\geq 3$
is odd and show
a relation of such codes and binary nonlinear Preparata codes.

\section{The space $X^n$ and the edges in $\{0,1\}^n$}\label{2}

An $edge$ in $\{0,1\}^n$ is a pair
$\{x,x'\}\subset \{0,1\}^n$ with $d_H(x,x')=1$. We say
that an edge $\{x,x'\}$ has {\it direction} $j$, iff
$x$ and $x'$ differ in the $j$th coordinate. Edges with identical
directions are called {\it parallel}.

We define a natural one-to-one mapping from the set of
edges in $\{0,1\}^n$ to $X^n$:
\begin{eqnarray*}
 \chi(\{(\sigma_1,\ldots,\sigma_{i-1},0,\sigma_{i+1},\ldots,\sigma_n),
 (\sigma_1,\ldots,\sigma_{i-1},1,\sigma_{i+1},\ldots,\sigma_n)\}) &&\\
\triangleq
 (\sigma_1,\ldots,\sigma_{i-1},*,\sigma_{i+1},\ldots,\sigma_n). &&
\end{eqnarray*}

We have the following straightforward proposition.
\bpro\label{dxde}
The distance $d_*(\cdot,\cdot)\triangleq d_H(\chi(\cdot),\chi(\cdot))$
between edges in $\{0,1\}^n$
 corresponding to the distance
$d_H$ in $X^n$ can be defined by the following rules.\\
 {\rm a)} $d_*(\{x,x'\},\{y,y'\})=
          \min(d_H(x,y),d_H(x,y'))$,
           if  $\{x,x'\}$ and $\{y,y'\}$ are parallel;\\
 {\rm b)} $d_*(\{x,x'\},\{y,y'\})=
          \min(d_H(x,y),d_H(x,y'))+2$
          otherwise.
\epro

In the figures that follow in the rest of the paper we use the edges interpretation of $X^n$.

A set $M$ of edges is called a {\it matching} iff the edges of $M$
are pairwise disjoint. A matching is {\it perfect} iff it covers
all $\{0,1\}^n$. All the codes in $X^n$ considered in Sections \ref{3}--\ref{5}
correspond to some matchings in $\{0,1\}^n$. The distance $3$
codes constructed in Section \ref{3} correspond  to perfect
matchings.

\section{Perfect distance $3$ codes}\label{3}

Let $  z\in \{0,1,*\}^n$ be an $n$-word with at most one symbol $*$.
Define
\begin{equation}\label{B}
B_{  z}\triangleq \{  y\in X^n|d_H(  y,  z)\leq 1\}.
\end{equation}

\begin{figure}
 \begin{center}
 \includegraphics{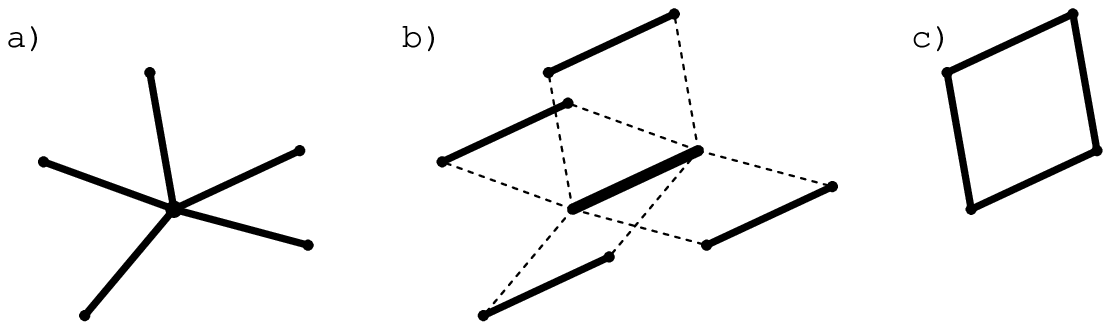}
 \caption{\label{f1} Diameter $2$ anticodes: a) $B_z$, $  z\in \{0,1\}^n$, b) $B_z$, $  z\in X^n$, c) square }
 \end{center}
\end{figure}

In
both cases  $  z\in \{0,1\}^n$ (see Figure \ref{f1},a)
and $  z\in X^n$  (see Figure \ref{f1},b) it is true that
\begin{itemize}
  \item $|B_{  z}|=n$;
  \item $B_{  z}$ is a diameter $2$ anticode.
\end{itemize}
Note that if $  z\in X^n$, then $B_{  z}$ is a radius $1$ ball in $X^n$.

\bexs Prove that if $n\geq 4$, then
$B_{z}$ is a maximal diameter $2$ anticode.
\eexs
\bnote The {\it square} diameter $2$ anticode
\[
\{(*,0,0,\ldots,0),(0,*,0,\ldots,0),(*,1,0,\ldots,0),(1,*,0,\ldots,0)\}
\]
{\rm (}see Figure \ref{f1},c{\rm )} is maximal for $n\leq 4$.
\enote

Assume $C$ is a diameter perfect distance $3$ code. Then
$|C|=|X^n|/|B_{  z}|=2^{n-1}$.
It is easy to see that in this case the balls $B_{  z}$, ${  z}\in C$ are
pairwise disjoint and cover all the space $X^n$.
This means that the distance $3$ code $C$ is {\it perfect} in the usual (not diameter) sense.
Moreover,
 if $C_i$ (respectively $X^n_i$) consists of all the words of $C$ (respectively $X^n$)
 with $*$ in the $i$th position,
then $\bigcup_{{  z}\in C_i}B_{  z}=X^n_i$. This means that $2^{n-1}=|X^n_i|$
is divisible by $|B_{  z}|=n$, i.\,e., $n=2^m$ for some $m\geq 2$.

In \cite{Svanstr:PhD,Svanstr:1999,vLinTol:1999} a construction
of a perfect distance $3$ constant-weight ternary code of length $n=2^m$ for each $m$ is
presented. The construction is based on the cosets of the cyclic Hamming binary code.

Other constructions presented in
\cite{
Kro:2001:Ternary}
are similar of some combinatorial constructions
of nonlinear perfect binary codes (with parameters of Hamming code) adapted for
perfect distance $3$ codes in $X^n$. Using these inductive constructions the
lower bound $2^{2^{n/2-2}}$ on the number of such codes is established. So, the
number of perfect codes in $X^n$ can be in some sense compared with the number
of perfect binary codes (see \cite{KroAvg:2008:lb} for a recent lower bound).

Subsection \ref{pd3pm} gives the interpretation of
perfect distance $3$ codes in $X^n$ in terms of edges in $\{0,1\}^n$.
In Subsection \ref{sscon} we present
a way to construct different perfect
distance $3$ codes in $X^n$ from cosets of Hamming code.
In the rest of the section we prove that this approach can give
nonequivalent codes.

\subsection{Perfect distance $3$ codes and perfect matchings in $\{0,1\}^n$}
\label{pd3pm}
Let
$\chi$ be the function defined in Section \ref{2}.
Then
\bpro\label{ppm}
The set $C\subset X^n$ is a distance $3$ perfect code if and
only if $M\triangleq\chi^{-1}(C)$ is a perfect matching in $E^n$
without parallel edges at distance $1$ or $2$.
\epro
\proofr
Only if. Assume $C$ is a distance $3$ perfect code. Code distance $3$
means that\\
a) edges of $M$ are pairwise disjoint, i.\,e., $M$ is matching;\\
b) if two different
edges of $M$ are parallel, then the distance between any element of
one edge and any element of the other is at least $3$.\\
The cardinality $2^{n-1}$ of perfect code $C$
means that the matching $M$ is perfect.

If. If two edges of the matching $M$ are not parallel, then
by Proposition \ref{dxde} the distance between corresponding words of $C$
is at least $1+2=3$.
If they are parallel, then the distance is at least $3$ by the condition.
The perfectness of $C$ is a corollary of the cardinality $2^{n-1}$ of
the perfect matching $M$.
\proofend

\subsection{Construction of perfect distance $3$ codes}
\label{sscon}
Let $F^n\triangleq\{0,1\}^n$ be the space of all binary $n$-words with the Hamming
distance. Assume $n=2^m\geq 4$ and $\alpha_1,\alpha_2,...,\alpha_n$ are all elements
of $F^m$. An {\it extended Hamming code} is the set
\[  H\triangleq\{(x_1,x_2,...,x_n)\in F^n\,|\,\sum_{i=1}^n x_i=0,\ \sum_{i=1}^n
x_i\alpha_i=0^m\}\]
of cardinality $2^n/2n$ and code distance $4$.
Note that varying the enumeration of the elements in $F^n$ we can get different
but equivalent extended Hamming codes.
If $\beta\in F^m$, then the sets
\begin{eqnarray*}
H^0_{\beta}&\triangleq&\{(x_1,x_2,...,x_n)\in F^n\,|\,\sum_{i=1}^n x_i=0,\ \sum_{i=1}^n x_i\alpha_i=\beta\},\\
H^1_{\beta}&\triangleq&\{(x_1,x_2,...,x_n)\in F^n\,|\,\sum_{i=1}^n x_i=1,\ \sum_{i=1}^n x_i\alpha_i=\beta\}
\end{eqnarray*}
are {\it even } and {\it odd cosets
} of $H$, respectively.
It follows from the cardinality of $H$ and the code distance $4$ that
each odd (even) word of $F^n$ ``sees'' exactly one word of
$H^0_{\beta}$ (respectively $H^1_{\beta}$) at the distance of $1$.
So, if $P$ and $Q$ are even and odd cosets of (may be different) Hamming codes,
then there is a unique function $q=q_{P,Q}:P\to Q$ such that $d_H(p,q(p))=1$ for
each $p\in P$. We also can define a function $r=r_{P,Q}:P\to X^n$ by the rule
$d_H(p,r(p))=d_H(q(p),r(p))=1$.
In terms of Section \ref{2} $r$ can be alternatively defined as
$r(p)\triangleq \chi(\{p,q(p)\})$.
Define $R(P,Q)\triangleq \{r_{P,Q}(p)\,|\,p\in P\}$.

\bpro\label{RPQ3} The code distance of $R(P,Q)$ is no less than $3$.
\epro
\proofr
Let $p$, $p'$ be different words from $P$. Then $d_H(p,p')\geq 4$.
Hence, $d_H(r(p),p')\geq 4$,
and $d_H(r(p),r(p'))\geq d_H(r(p),p')-d_H(p',r(p'))\geq 3$.
\proofend

Let $f:F^m\to F^m$ be a linear operator with the following properties:\\
a) $f$ is a one-to-one operator,\\
b) $f+Id$ is a one-to-one operator.
\btheorem\label{th1}
The set $C_{H,f}\triangleq\bigcup_{\beta\in F^m}R(H^0_{\beta},H^1_{f(\beta)})$ is a perfect
distance $3$ code in $X^n$.
\etheorem
\proofr
Consider the set of edges
\[
M\triangleq\chi^{-1}(C_{H,f})=\bigcup_{\beta\in F^m} M_{\beta},
\quad \mbox{where } M_\beta \triangleq \chi^{-1}(R(H^0_{\beta},H^1_{f(\beta)})).
\]

Each edge in $M_\beta$ consists of a word in $H^0_{\beta}$ and
 a word in $H^1_{f(\beta)}$.
Since $f$ is one-to-one, the edges in $M$ are pairwise disjoint.
Since $|M|=2^m|H|=2^{n-1}$, $M$ is a perfect matching.
For fixed $\beta$ all edges
in $M_\beta$ have the same direction $i_{\beta}$
defined by equality $\alpha_{i_\beta}=\beta+f(\beta)$.
Since $f+Id$ is one-to-one, the directions $i_\beta$, ${\beta\in F^m}$,
are pairwise different. So, if two different edges in $M$ are parallel, then
they belong to a same set $M_\beta$ and the distance between them is
at least $3$ by Proposition \ref{RPQ3}.
By Proposition \ref{ppm} the theorem is proved.
\proofend

\bexam The linear operators defined by the matrices
\[
  \left\lgroup
  {1 \ 1 \atop 1 \ 0  }
  \right\rgroup
  ,\  m=2 \quad \mbox{ and }  \quad
  \left\lgroup
  \begin{array}{ccc}
    1 & 1 & 0 \\
    0 & 1 & 1 \\
    1 & 0 & 0
  \end{array}
  \right\rgroup
,\  m=3
\]
satisfy conditions {\rm a)} and {\rm b)}.
\eexam
\bexam\label{ex2}
Let $f':F^{m'}\to F^{m'}$ and $f'':F^{m''}\to F^{m''}$ be linear operators
satisfying conditions {\rm a)} and {\rm b)}.
Then the operator $g:F^{m'+m''}\to F^{m'+m''}$ defined by
$g(x',x'') \triangleq (f'(x'),f''(x''))$ also satisfies {\rm a)} and {\rm b)}.
\eexam
\bexam\label{ex3}
Let us consider words of $F^n$ as a representation of elements of $GF(2^m)$.
Let $\gamma$ be a primitive element of $GF(2^m)$.
Then the operator $h:GF(2^m)\to GF(2^m)$ defined by
$h(x)=\gamma x$ satisfies {\rm a)} and {\rm b)}.
\eexam
In Subsection \ref{ssnneq} we prove that the codes $C_{H,g}$ and
$C_{H,h}$ corresponding to the operators $g$ and $h$ from Examples
\ref{ex2} and \ref{ex3} cannot be equivalent.

\subsection{Automorphism group and transitivity}
Let $\pi:\{1,\ldots,n\}\leftrightarrow\{1,\ldots,n\}$
be a coordinate permutation and $z\in\{0,1\}^n$.
A pair $\tau=(\pi,z)$ is called an {\it automorphism} of a code
$C\subset X^n$ iff $x\in C$ implies
$\tau(x)\in C$, where $\tau(x)\triangleq\pi(x)+z$,
$\pi(x_1,\ldots,x_n)\triangleq (x_{\pi^{-1}(1)},\ldots,x_{\pi^{-1}(n)})$,
and $+$ is a coordinate-wise mod~$2$ addition extended by equalities
$*+1=*+0=*$.
It is clear that if $\tau=(\pi,z)$ is an automorphism of $C$, then
$\tau^{-1}=(\pi^{-1},\pi^{-1}(z))$ is also an automorphism of $C$.
Denote the group of all automorphisms of a
code $C$ by $\Aut(C)$. We say that $\Aut(C)$ is {\it transitive} iff for each
$x,y\in C$ there is $\tau\in \Aut(C)$ such that $\tau(x)=y$.

\btheorem\label{ptra}
The automorphism group of the code $C_{H,f}$
defined in Theorem \ref{th1} is transitive.
\etheorem
\proofr
1) First we will show that for each even weight vector $z\in \{0,1\}^n$
there is $\tau_z\in\Aut(C_{H,f})$ which sends $0^n$ to $z$.
Define $\tau_z\triangleq \pi_z+z$,
where $\pi_z$ is defined by identity
\begin{equation}\label{dpi}
\alpha_{\pi_z(i)}=\alpha_i+s_z+f(s_z),
\quad s_z\triangleq \sum_{j=1}^n z_j\alpha_j,
\end{equation}
i.\,e., $\pi_z(i)$ is the number of
$\alpha_i+s_z+f(s_z)$
in the numeration $\{\alpha_1,\ldots,\alpha_n\}$ of elements of
$F^m$. It is obvious that $\tau_z(0^n)=z$, and we only need to
check that $\tau_z$ is an automorphism of $C_{H,f}$.

Let $r=\tau_z(r')$ and $r\in C_{H,f}$,
i.\,e, $r\in R(H^0_{\beta},H^1_{f(\beta)})$ for some
$\beta\in F^m$.
Let $p\triangleq r^{-1}_{H^0_{\beta},H^1_{f(\beta)}}(r)$,
and $q\triangleq q_{H^0_{\beta},H^1_{f(\beta)}}(p)$.
Then $p\in H^0_{\beta}$ and $q\in H^1_{f(\beta)}$.
If we show that $p'\triangleq\tau_z^{-1}(p)\in H^0_{\beta'}$
and $q'\triangleq\tau_z^{-1}(q)\in H^1_{\beta'}$
for some $\beta'\in F^m$,
then we get $q'=q_{H^0_{\beta'},H^1_{f(\beta')}}(p')$
and $r'=r_{H^0_{\beta'},H^1_{f(\beta')}}(p')\in C_{H,f}$,
which implies that the statement of the theorem is true.

The condition $p\in H^0_{\beta}$ means that
\[ \sum_{i=1}^n{p_i\alpha_i}=\beta,\quad \sum_{i=1}^n{p_i}=0. \]
Since $p=\tau_z(p')=\pi_z(p')+z$ and $z$ is even weight vector we have
\begin{eqnarray}
 \sum_{i=1}^n{p'_{\pi_z^{-1}(i)}\alpha_i}=\sum_{i=1}^n{z_i\alpha_i}+\beta,\quad \sum_{i=1}^n{p'_i}=0,
 \nonumber \\
 \sum_{i=1}^n{p'_{i}\alpha_{\pi_z(i)}}=\sum_{i=1}^n{z_i\alpha_i}+\beta,\quad \sum_{i=1}^n{p'_i}=0.
\nonumber\end{eqnarray}
Substituting (\ref{dpi}) we get
\begin{eqnarray}
 \sum_{i=1}^n{p'_{i}\alpha_{i}}+\sum_{i=1}^n{p'_{i}(s_z+f(s_z))}
=s_z+\beta,\quad \sum_{i=1}^n{p'_i}=0, \nonumber \\
 \sum_{i=1}^n{p'_{i}\alpha_{i}}
=s_z+\beta,\quad \sum_{i=1}^n{p'_i}=0,
\nonumber\end{eqnarray}
i.\,e., $p'\in H^0_{\beta'}$ for
$\beta'\triangleq s_z+\beta$.

The condition $q\in H^1_{f(\beta)}$ means that
\[ \sum_{i=1}^n{q_i\alpha_i}=f(\beta),\quad \sum_{i=1}^n{q_i}=1. \]
Since $q=\tau_z(q')=\pi_z(q')+z$ and $z$ is an even weight vector we have
\begin{eqnarray}
 \sum_{i=1}^n{q'_{\pi_z^{-1}(i)}\alpha_i}=\sum_{i=1}^n{z_i\alpha_i}+f(\beta),\quad \sum_{i=1}^n{q'_i}=1,
 \nonumber \\
 \sum_{i=1}^n{q'_{i}\alpha_{\pi_z(i)}}=\sum_{i=1}^n{z_i\alpha_i}+f(\beta),\quad \sum_{i=1}^n{q'_i}=1.
\nonumber\end{eqnarray}
Substituting (\ref{dpi}) we get
\begin{eqnarray}
 \sum_{i=1}^n{q'_{i}\alpha_{i}}+\sum_{i=1}^n{q'_{i}(s_z+f(s_z))}
=s_z+f(\beta),\quad \sum_{i=1}^n{q'_i}=1, \nonumber\\
 \sum_{i=1}^n{q'_{i}\alpha_{i}}
=(s_z+f(s_z))
+s_z+f(\beta)=f(s_z+\beta),\quad \sum_{i=1}^n{q'_i}=1,
\nonumber\end{eqnarray}
i.\,e., $q'\in H^1_{f(\beta')}$.

The condition $d_H(p,q)=1$ means $d_H(p',q')=1$. So, we get
$q'=q_{H^0_{\beta'},H^1_{f(\beta')}}(p')$ and, as a corollary,
$r'=r_{H^0_{\beta'},H^1_{f(\beta')}}(p')
\in R(H^0_{\beta'},H^1_{f(\beta')})\subset C_{H,f}$.
Therefore $\tau_z\in \Aut(C_{H,f})$.

2) Let $r,r'$ be arbitrary words in $C_{H,f}$ and $r=r(p)$, $r'=r(p')$
where $p,p'\in F^n$ are even weight vectors.
Define $\tau=\tau_{p'}\tau_p^{-1}\in \Aut(C_{H,f})$.
Then $\tau(p)=p'$. Since there is only one element in $C_{H,f}$ at
the distance $1$ from $p'$, we also have $\tau(r)=r'$. This proves
the transitivity of $\Aut(C_{H,f})$.
\proofend

\subsection{Nonequivalence}\label{ssnneq}
Two codes $C$ and $C'$ in $X^n$ are {\it equivalent}
iff $C=\pi(C')+z$
for some coordinate permutation $\pi$ and $z\in F^n$.
The goal of this subsection is to show that the construction
of Subsection \ref{sscon} can give nonequivalent distance
$3$ perfect codes in $X^n$.

Let, as in Example \ref{ex2}, $f':F^{m'}\to F^{m'}$
and $f'':F^{m''}\to F^{m''}$ be linear operators
satisfying conditions a) and b),
and the operator $g:F^{m'+m''}\to F^{m'+m''}$ be defined by
$g(x',x'')\triangleq (f'(x'),f''(x''))$. Let $n'=2^{m'}$,
$n''=2^{m''}$, and
$n=2^m$, where $m=m'+m''$.

\bpro\label{pro+} Fixing $n-n'$ coordinates by zeroes we can get
a perfect distance $3$ code in $X^{n'}$ from the
perfect distance $3$ code $C_{H,g}\subset X^n$.
\epro
\proofr
W.l.o.g. assume that $\alpha_1,\ldots,\alpha_{n'}$ have zeroes
in the last $m''$ positions,
i.\,e., $\{\alpha_1,\ldots,\alpha_{n'}\}=F^{m'}\times \{0\}^{m''}$.

It is enough to show that there are $2^{n'-1}$
codewords in $C_{H,g}$ with zeroes in the last $n-n'$ coordinates.

Choose an arbitrary even weight vector $p\in F^{n'}\times\{0\}^{n-n'}$.
Let $\beta\triangleq\sum_{i=1}^n p_i\alpha_i=\sum_{i=1}^{n'}
p_i\alpha_i$.
Then $\beta\in F^{m'}\times \{0\}^{m''}$
and, by definition of $g$,
we also have $g(\beta)\in F^{m'}\times \{0\}^{m''}$.
Let $q\triangleq q_{H^0_\beta,H^1_{g(\beta)}}(p)$.
Assume $q\not\in F^{n'}\times\{0\}^{n-n'}$.
Since $d_H(p,q)=1$, the word $q$ has exactly one nonzero
coordinate $j$ larger than $n'$.
But this contradicts the equation
$g(\beta)=\sum_{i=1}^n q_i\alpha_i=\sum_{i=1}^{n'} q_i\alpha_i+\alpha_j$
because $\alpha_j\not\in F^{m'}\times \{0\}^{m''}$.
So, $q\in F^{n'}\times\{0\}^{n-n'}$ and, consequently,
$r\triangleq r_{H^0_\beta,H^1_{g(\beta)}}(p)\in
X^{n'}\times\{0\}^{n-n'}$. Since there are $2^{n'-1}$ ways to
choose $p$, we get $2^{n'-1}$ codewords $r\triangleq r_{H^0_\beta,H^1_{g(\beta)}}(p)\in
X^{n'}\times\{0\}^{n-n'}\cap C_{H,g}$. After deleting the last
$n-n'$ zeroes in these codewords we get a distance $3$ code in $X^{n'}$ of
cardinality $2^{n'-1}$, i.\,e., a perfect code.
\proofend

As in Example \ref{ex3},
consider the words of $F^n$ as representation of the elements of $GF(2^m)$.
Define an operator $h:GF(2^m)\to GF(2^m)$ by
$h(x)\triangleq\gamma x$, where $\gamma$ is a primitive element of $GF(2^m)$.

\bpro\label{pro-} Fixing $n-k$ coordinates, $1<k<n$, we cannot get
a perfect distance $3$ code in $X^{k}$ from the
perfect distance $3$ code $C_{H,h}\subset X^n$.
\epro
\proofr
Denote $M=\chi^{-1}(C_{H,h})$.
Assume the statement is not true. This means that $M$ contains
(as a subset) a matching $M^k$ from $2^{k-1}$ edges which covers all the
words of $F^n$ with fixed $n-k$ coordinates.

W.l.o.g. we assume $\alpha_1=0^m$.

1) First we consider the case when $M^k$ contains an edge of the first direction.
Note that $M^k$ also contains
an edge of another direction, otherwise $k=1$.
Let this edge be $\{a^0, a^1\}$. Further we construct
a sequence $a^0,a^1,a^2,\ldots$ by the following inductive
rules:\\
a) $a^{2j+2}= a^{2j+1}+(1,0,\ldots,0)$, $j=0,1,\ldots$;\\
b) $\{a^{2j},a^{2j+1}\}\in M$, $j=0,1,\ldots$
(this rule is correct because $M$
is a perfect matching by Proposition \ref{ppm}).\\
It is easy to see that by induction
$\{a^{2j},a^{2j+1}\}\in M^k$, $j=0,1,\ldots$.
Let $i_j$ be the direction of the edge $\{a^{2j},a^{2j+1}\}$.
Then $\alpha_{i_j}=s_{a^{2j}}+s_{a^{2j+1}}$ where
$s_{x}\triangleq\sum_{i=1}^{n}x_i\alpha_i$.
By the construction $s_{a^{2j+1}}=\gamma s_{a^{2j}}$ and
$s_{a^{2j+2}}= s_{a^{2j+1}}$. Therefore
$\alpha_{i_{j+1}}=\gamma\alpha_{i_j}$ and, by induction,
$\alpha_{i_{j}}=\gamma^j\alpha_{i_0}$.
Since $\gamma$ is a primitive element of $GF(2^m)$,
we have
$\{\alpha_{i_0},\ldots,\alpha_{i_{n-2}}\}=\{\alpha_2,\ldots,\alpha_n\}$
(recall that $\alpha_1$ is the zero element).
Consequently, the directions of $\{a^{2i},a^{2i+1}\}$,
$i=0,\ldots,n-2$, are pairwise different, and $M^k$ contains edges
of all directions. This contradicts the condition $k<n$.

2) Assume $M^k$ contains an edge $e$ of the $l$th direction
for some $l\in\{1,\ldots,n\}$.
By Theorem \ref{ptra} there is $\tau\in\Aut(C_{H,h})$
such that $\tau(e)$ has the first direction. Then $\tau(M^k)\subset M$ is
again a matching from $2^{k-1}$ edges which covers all the
words of $F^n$ with fixed $n-k$ coordinates. And we get a
contradiction by p.1.
\proofend

As obviously follows from Propositions \ref{pro+} and \ref{pro-},
\bcorol
The codes $C_{H,g}$ and $C_{H,h}$ are nonequivalent.
\ecorol
Thereby we have constructed nonequivalent perfect distance $3$
codes in $X^n$ based on the cosets of the same Hamming code.

\section{Diameter perfect distance $4$ codes}\label{4}
\begin{figure}
 \begin{center}
 \includegraphics{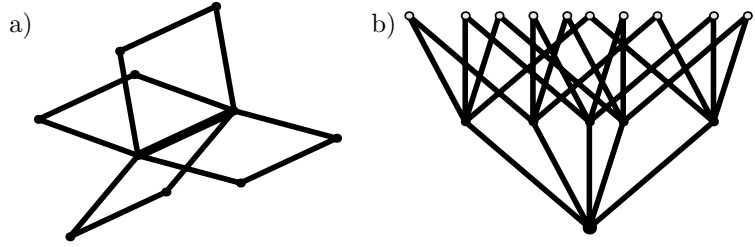}
 \caption{\label{f2} a) Diameter $3$ anticode; b) diameter $4$ anticode}
 \end{center}
\end{figure}

Let $z\in X^n$ and $z_0$, $z_1$ be different words in $F^n$ at the distance
$1$ from $z$. Define $A_z\triangleq B_{z}\cup B_{z_0}\cup B_{z_1}$ (see Figure \ref{f2},a)
where $B_{z}, B_{z_0}, B_{z_1}$ are
defined by (\ref{B}). Then
\begin{itemize}
  \item $|A_{ z}|=3n-2$;
  \item $A_{ z}$ is a diameter $3$ anticode.
\end{itemize}
\bexs Prove that for a sufficiently large $n$ the set $A_{z}$
is a maximal diameter $3$ anticode.
\eexs
So, $n2^{n-1}\vdots (3n-2)$, i.\,e., $n=(2^{2m}+2)/3$ is a necessary condition
for existing diameter perfect distance $4$ codes.

For $m=2$ such a code of cardinality $12$ exists  \cite{Svanstr:PhD} and can
be constructed from the rows of the $6\times 6$ conference matrix
(see, e.\,g., \cite{MWS} for the definition)
\[
\left\lgroup
  \begin{array}{cccccc}
    0&\hphantom{-}1&\hphantom{-}1&\hphantom{-}1&\hphantom{-}1&\hphantom{-}1\\
    1&\hphantom{-}0&\hphantom{-}1&-1&-1&\hphantom{-}1\\
    1&\hphantom{-}1&\hphantom{-}0&\hphantom{-}1&-1&-1\\
    1&-1&\hphantom{-}1&\hphantom{-}0&\hphantom{-}1&-1\\
    1&-1&-1&\hphantom{-}1&\hphantom{-}0&\hphantom{-}1\\
    1&\hphantom{-}1&-1&-1&\hphantom{-}1&\hphantom{-}0
  \end{array}
\right\rgroup
\]
by substitutions $0\to *$, $1\to 0$, $-1\to 1$
and $0\to *$, $1\to 1$, $-1\to 0$.


Assume $C$ is a diameter perfect distance $4$ code in $X^n$.
Then $|C|=|X^n|/|A_{z}|=\frac{n2^{n-1}}{3n-2}$.
Define $C_1=\{x\in F^n \,|\,d_H(x,C)=1\}$ and $C_2=F^n \backslash C_1$.
It is easy to see that every word from $C_1$ has exactly
$1$ neighbor from $C_1$ and $n-1$ neighbors from $C_2$;
every word from $C_2$ has at most $n/2$ neighbors from $C_1$.
Since $|C_1|/|C_2|=2|C|/(2^n-2|C|)=\frac{n}{2(n-1)}$,
every word from $C_2$ has exactly $n/2$ neighbors from $C_1$
and $n/2$ neighbors from $C_2$.
This means that $\{C_1,C_2\}$ is a so-called \emph{perfect coloring}, 
or \emph{equitable partition}, of $F^n$
with the parameters $((1,n-1),(n/2,n/2))$.
Unfortunately, such perfect colorings exist only
for $n=6$ \cite{FDF:CorrImmBound}.
So, this value is the only one for which a diameter perfect distance $4$ code in $X^n$ exists.

\section{Diameter perfect distance $5$ codes}\label{5}

An operator $f:F^m\to F^m$ is {\it APN} ({\it almost perfect nonlinear})
iff the system of equations
\[ \cases{x + y = a \cr f(x) + f(y) = b} \]
has $0$ or $2$ solutions $(x; y)$ for every $(a; b)\neq (0^m; 0^m)$.

\begin{lmm}[see \cite{CarChaZin}]
If $f:F^m\to F^m$ is an {\it APN} operator, then the binary code $C_f\subset F^n$ $(n={2^m})$ defined by
\begin{equation}\label{BCH}
  C_f\triangleq\{(x_1,x_2,...,x_n)\in F^n\,|\,\sum_{i=1}^n x_i=0,\ \sum_{i=1}^n x_i\alpha_i=0^m,
\ \sum_{i=1}^n x_if(\alpha_i)=0^m\}
\end{equation}
has a code distance $6$.
\end{lmm}
Recall that $\alpha_1,\ldots,\alpha_n$ are all the elements of $F^m$.
Assume an APN operator $f:F^m\to F^m$ is one-to-one. Then the sets $P$ and $P'$
defined as
\begin{eqnarray*}
  P&\triangleq&\{(x_1,x_2,...,x_n)\in F^n\,|\,\sum_{i=1}^n x_i=0,\ \sum_{i=1}^n
  x_i\alpha_i=0^m\},
\\
  P'&\triangleq&\{(x_1,x_2,...,x_n)\in F^n\,|\,\sum_{i=1}^n x_i=0,\ \sum_{i=1}^n x_i
  f(\alpha_i)=0^m\}
\end{eqnarray*}
are two different Hamming codes, and their intersection is the distance $6$
code $C_f$, see (\ref{BCH}). Let $Q$ be an odd coset of $P'$. Then
\btheorem
The set $R(P,Q)$ defined as in Section \ref{3} is a diameter perfect distance
$5$ code in $X^n$.
\etheorem
\proofr Let $r(x)$ and $r(x')$ be different words in
$R(P,Q)$, where $x,x'\in P$ and $d_H(x,r(x))=d_H(x',r(x'))=1$.
(The functions $q=q_{P,Q}$ and $r=r_{P,Q}$ are defined in Section \ref{3}.)

1. Assume $r(x)$ and $r(x')$ contain $*$ in a same position. It can be shown that $x$ and $x'$ are
in a same coset of $C_f$ in this case. Hence, $d_H(r(x),x')\geq d_H(x,x')\geq 6$ and
$d_H(r(x),r(x'))\geq d_H(r(x),x')-d_H(x',r(x'))\geq 6-1=5$.

2. Assume $r(x)$ and $r(x')$ contain $*$ in different positions. If
$d_H(x,x')\geq 6$, then $d_H(r(x),r(x'))\geq 5$, see p. 1.
Let $d_H(x,x')=4$ and $x$, $x'$ differ in positions $i_1,i_2,i_3,i_4$.

2.1. If both $r(x)$ and $r(x')$ contain $*$ in positions from $\{i_1,i_2,i_3,i_4\}$
(for example, $i_1,i_2$),
then $q(x)$ and $q(x')$ differ in only two positions (respectively $i_3,i_4$).
This contradicts the fact that $q(x),q(x')\in Q$.

2.2. If at least one word of $r(x),r(x')$ contains $*$ in position $j$ different
from $i_1,i_2,i_3,i_4$, then $r(x),r(x')$ differ in at least five positions
$j,i_1,i_2,i_3,i_4$.

It remains to prove that the distance $5$ code $R(P,Q)$ is diameter
perfect. It is enough to show the existence of a diameter $4$ anticode of
cardinality $|X^n|\cdot |R(P,Q)|^{-1}=n2^{n-1}\cdot (2^n/2n)^{-1}=n^2$ in $X^n$.
The set $\{y\in X^n|d_H(0^n,y)\leq 2\}$ (see Figure \ref{f2},b)
satisfies these conditions.
\proofend

It is known that if $m$ is odd, then APN one-to-one operators exist, see the
following examples. So, we get a construction of
diameter perfect distance $5$ codes in $X^{2^m}$ for odd $m$.

\bnote
The optimal distance $5$ code in $X^8$ belongs to the series of distance $(p^m+3)/2$
codes in $X^{p^m+1}$ constructed from Jacobsthal matrices {\rm \cite{Svanstr:PhD,OstSva2002}}.
\enote

\bexam
The operator $u:GF(2^m)\to GF(2^m)$ defined by $u(\alpha)=\alpha^3$
for any $\alpha\in GF(2^m)$ is APN.
But $u$ is one-to-one if and only if $m$ is odd.
\eexam

\bexam
The operator $v:GF(2^m)\to GF(2^m)$ defined by
\begin{eqnarray*}
v(\alpha)&=&\alpha^{-1},\quad \alpha\neq 0, \\
v(0)&=&0
\end{eqnarray*}
is one-to-one.
But it is APN if and only if $m$ is odd.
\eexam

The corresponding binary codes $C_u$ and $C_v$ are distance $6$ BCH and reversible Melas
codes (see \cite{MWS}).

{\bf Question.} Do diameter perfect distance $5$ codes in $X^{2^m}$
exist for even $m$?

\subsection{Connection with Preparata codes}\label{sspre}
Using the notation of Section \ref{5} 
we will construct binary nonlinear codes
with parameters of Preparata codes (see \cite{MWS}) of length
$2n$, cardinality $2^{2n}/4n^2$ and distance $6$. Our construction
is in one step from the representation of such codes given in
\cite{BvLW}, see also \cite{vDamFDFla:2003}.
The goal is to illustrate connection
between constructions of distance $5$ diameter perfect
constant-weight ternary codes and optimal binary Preparata-like codes.

Assume $f:F^m\to F^m$ is a one-to-one APN operator which satisfies the
following additional property:
\begin{eqnarray}\label{propf}
&&\mbox{the system}\quad
\cases{a+b+c+d=0 \cr a+b+e+t=0 \cr
f(a)+f(b)+f(c)+f(d)+f(e)+f(t)=0}
\\ \nonumber
&&\mbox{has no solutions with pairwise different $a,b,c,d,e,t\in F^m$.}
\end{eqnarray}

Let, as in Section \ref{5},
\begin{eqnarray*}
  P&\triangleq&\{(x_1,x_2,...,x_n)\in F^n\,|\,\sum_{i=1}^n x_i=0,\ \sum_{i=1}^n
  x_i\alpha_i=0^m\},
\\
  Q&\triangleq&\{(x_1,x_2,...,x_n)\in F^n\,|\,\sum_{i=1}^n x_i=1,\ \sum_{i=1}^n x_i
  f(\alpha_i)=0^m\},
\end{eqnarray*}
and the operator $q=q_{P,Q}:P\to Q$ be defined by
$\forall p\in P: d_H(p,q(p))=1$ as in Section \ref{3}.

\btheorem\label{thpr}
The code
\begin{equation}\label{prep}
\overline{P}\triangleq\{ (a+b+q(b),a+q(b))\,:\,a,b\in P \}
\end{equation}
is a binary distance $6$ code of length $2n$ and cardinality $2^{2n-2}/n^2$.
\etheorem
\proofr It is obvious that $|\overline{P}|=|P|^2=2^{2n-2}/n^2$.
Let $c=(a+b+q(b),a+q(b))$ and $c'=(a'+b'+q(b'),a'+q(b'))$
be different codewords in $\overline{P}$, where $a,b,a',b'\in P$.
\def\dc{\tilde c}
\def\da{\tilde a }
\def\db{\tilde b }
\def\dq{\tilde q }
\def\dv{\tilde \omega }

Denote $\dc\triangleq c+c'$, $\da\triangleq a+a'$, $\db\triangleq b+b'$,
$\dq\triangleq q(b)+q(b')$, $\dv\triangleq (b+q(b))+(b'+q(b'))$,
$w_H(x)\triangleq d_H(x,(0\ldots 0))$.
Note that\\
a) $\da,\db \in P$,\\
b) $\dq\in P'$,\\
c) $\dv$ has $0$ or $2$ nonzero coordinates,\\
d) $\dv=0^n$ iff $\db=\dq\in C_f$, where $C_f$ is defined in (\ref{BCH}),\\
e) $\db=0^n$ iff $\dq=0^n$.

1) Let $w_H(\db)=0$. Then $w_H(\da)\geq 4$ and
$d_H(c,c')=w_H(\dc)=w_H(\da,\da)\geq 4+4=8$, see (a,e).

2) Let $w_H(\db)\geq 6$.
Then $w_H(\dc)=w_H(\da+\db+\dq,\da+\dq)
\geq w_H((\da+\db+\dq)+(\da+\dq))=w_H(\db)\geq 6$.

3) Let $w_H(\da)=0$.
Then $w_H(\dc)=w_H(\dv,\dq)\geq min(2+4,0+6)=6$, see (b,c,d,e).

4) Let $w_H(\da)\geq 6$.
If $\da\neq \dq$,
then $w_H(\dc)=w_H(\da+\dv,\da+\dq)\geq (6-2)+2=6$, see (c).
Otherwise $\dq=\da\in P$.
Since $\db\in P$ and $d_H(\db,\dq)=w_H(\dv)\leq 2$, we get $\db=\dq=\da$,
and $w_H(\dc)\geq 6$ by p.2).

5) Let $w_H(\db)=4$ and $w_H(\da)=4$, the last case. Consider
subcases taking into account that $\da+\db\in P$.

5.1) Assume $w_H(\da+\db)\geq 6$. Then
$w_H(\dc)=w_H(\da+\dv,\da+\db+\dv)\geq (4-2)+(6-2) = 6$, see (c).

5.2) Assume $\da=\db$. Then
$w_H(\dc)=w_H(\dq,\dv)\geq 4+2 = 6$, see (b,c,d,e).

5.3) Assume $w_H(\da+\db)=4$.
We have $\da,\db,\da+\db\in P$ ,$\dq\in P'$
and $P\cap P'$ is the code $C_f$ with distance $6$.
Consequently, $d_H(\da,\dq)\geq 2$,
$d_H(\da+\db,\dq)\geq 2$, and $d_H(\db,\dq)= 2$.
If $d_H(\da,\dq)> 2$ or $d_H(\da+\db,\dq)> 2$, then
$w_H(\dc)=w_H(\da+\db+\dq,\da+\dq)\geq 6$.
If $d_H(\da,\dq)=d_H(\da+\db,\dq)=d_H(\db,\dq)= 2$, then
the only possible case is $w_H(\dq)=6$ and $\dq$ has ones in
the positions in which $\da$ or $\db$ has one. But by the construction
of $P$ and $P'$ and property (\ref{propf}) this is impossible.

So, the code distance of $\overline{P}$ is no less than $6$.
Since $w_H(\da)=0$, $w_H(\dq)=4$ is a possible case, $6$ is a tight value.
\proofend

\bnote
Let $\omega(b)\triangleq b+q(b)$. Then {\rm (\ref{prep})} can be represented
as
\begin{equation}\label{eVF}
\overline{P}\triangleq\{ (a+\omega(b),a+b+\omega(b))\,:\,a,b\in P \}.
\end{equation}
It appears similar to the well known $(u,u+v)$-construction {\rm (}see {\rm \cite{MWS}}{\rm )}
improved by the weight
$1$ vector $\omega(b)$ defined as starting in $b\in P$ and finishing in $Q$.
\enote

\bpro\label{ex99}
Let $m$ be odd. Let $f:GF(2^m)\to GF(2^m)$ be defined by
$f(\alpha)=\alpha^{2^l+1}$ for any $\alpha \in GF(2^m)$
where $l$ and $m$ are relatively prime.
Then $f$ is a one-to-one APN operator satisfying {\rm (\ref{propf})}.
\epro
\proofr
We will show only (\ref{propf}) because the other is well known
and simple to prove.
We will use well known relations
\[ (a+b)^{2^l}=a^{2^l}+b^{2^l},\quad
(a+b)^{2^l+1}=a^{2^l+1}+a^{2^l}b+ab^{2^l}+b^{2^l+1}
\]
and
\begin{equation}\label{ne0}
   a^{2^l}+a+1\neq 0.
\end{equation}

Let the system in (\ref{propf})
has a solution with pairwise different $a,b,c,d,e,t\in F^m$.
W.l.o.g. we can assume $a=0\in GF(2^m)$, otherwise substituting
$a'=a+a=0$, $b'=b+a$, $c'=c+a$, $d'=d+a$, $e'=e$, $t'=t$ we get
that the system in (\ref{propf}) holds for $a'=0,b',c',d',e',t'$.
Further w.l.o.g. we assume $b=1$, otherwise we can consider
$a''=0/b=0$, $b''=b/b=1$,
$c''=c/b$, $d''=d/b$, $e''=e/b$, $t''=t/b$.

So, we have
\[ 1+c+d=0,\quad 1+e+t=0,\quad 1+c^{2^l+1}+d^{2^l+1}+e^{2^l+1}+t^{2^l+1}=0.
\]

Substituting $d=c+1$ and $t=e+1$ in the last equation we get
\[ (c+e)^{2^l}+(c+e)+1=0 \]
which contradicts (\ref{ne0}).
\proofend

Codes constructed in Theorem \ref{thpr} with $f$ from Proposition \ref{ex99}
are equivalent to Preparata-like codes from \cite{BvLW,Dumer:1976}.
The case $l=1$ corresponds to the original Preparata codes \cite{Preparata:1968}.
Unfortunately the operators from Proposition
\ref{ex99} fulfill the set of all known operators appropriate for
constructing Preparata-like codes by (\ref{prep}) or (\ref{eVF}).
Another known class of Preparata-like codes is the class of $Z_4$-linear Preparata
codes;
\cite{Tokareva:2005} gives a representation of such codes by formulas similar to (\ref{eVF}).

\bibliographystyle{plain}
\bibliography{../../../k}
\end{document}